\documentclass{gen-p-l}
\newcommand{\iso}{{\widetilde\to}}

\newcommand\BOX[1]{\mathop{\mathrm{#1}}}
\newcommand\BOXB[1]{\mathop{\mathbf{#1}}}

\newcommand{\supp}{\BOX{supp}}

\newcommand{\Hom}{\BOX{Hom}}
\newcommand{\End}{\BOX{End}}

\newcommand\C{\mathbb{C}}
\newcommand\Z{\mathbb{Z}}
\newcommand\R{\mathbb{R}}

\numberwithin{equation}{section}
\newtheorem{Th}{Theorem}[section]
\newtheorem{Pp}[Th]{Proposition}
\newtheorem{Lm}[Th]{Lemma}
\newtheorem{Co}[Th]{Corollary}

\theoremstyle{definition}
\newtheorem{Df}[Th]{Definition}
\newtheorem{Exm}[Th]{Example}

\theoremstyle{remark}
\newtheorem{Rem}{Remark}
\newcommand{\select}[1]{\textit{#1}}

\newcommand{\cal}{\mathcal}

\newcommand{\cF}{{\cal F}}
\newcommand{\cL}{{\cal L}}
\newcommand{\cP}{{\cal P}}
\newcommand{\cS}{{\cal S}}

\newcommand{\tF}{{\widetilde F}}
\newcommand{\tL}{{\widetilde L}}

\newcommand\dl{^{\vee}}
\newcommand\Sky{\BOXB{Sky}}
\newcommand\tSky{\widetilde{\Sky}}
\newcommand\LU{\BOXB{Loc}_u}

\newcommand\DR{\BOXB{DR}\nolimits}

\newcommand\Vect{\Omega^{-1}}
\newcommand\Diff{\BOX{Diff}}

\newcommand\Four{\BOXB{Four}\nolimits}

\newcommand\curv{\BOX{curv}}
\newcommand\mon{\BOX{mon}}
\newcommand\al{\bullet}
\newcommand\ex[2]{\BOX{Exp}(-2\pi\sqrt{-1}\langle #2,#1\rangle)}
\newcommand\toup[1]{\mathop{\to}\limits^{#1}} %???
\newcommand\tcL{\widetilde\cL}

\begin{document}

\title{Fukaya category and Fourier transform}
\author{D.~Arinkin}
\address{Department of Mathematics, Harvard University, Cambridge, Massachusetts 02138}
\email{arinkin@math.harvard.edu}
\author{A.~Polishchuk}
\address{Department of Mathematics, Harvard University, Cambridge, Massachusetts 02138}
\email{apolish@math.harvard.edu}

\subjclass{Primary 53C15; Secondary 14H52}
\date{November 4, 1998}%???

\begin{abstract}
We construct a version of Fourier transform for families of real tori. This 
transform defines a functor from certain category
associated with a symplectic family of tori to the category of holomorphic 
vector bundles on the dual family (the dual family has a natural complex
structure). In the $1$-dimensional case, the former category is closely 
related to the Fukaya category. 
\end{abstract}

\maketitle

\section*{Introduction}

This paper grew out from attempts to understand better the
homological mirror symmetry for elliptic curves.
The general homological mirror conjecture formulated by
M.~Kontsevich in \cite{Kon} asserts that the derived category
of coherent sheaves on a complex variety is equivalent to (the
derived category of) the Fukaya category of the mirror dual symplectic
manifold.
This equivalence was proved in \cite{PZ} for the case of elliptic curves
and dual symplectic tori.
However, the proof presented in \cite{PZ} is rather computational and does not
give a conceptual construction of a functor between two categories.
In the present paper we fill up this gap by providing such a
construction. We also get a glimpse of what is going on in the
higher-dimensional case.

The idea is to use a version of the Fourier transform for families of real tori
which generalizes the well-known correspondence between smooth functions
on a circle and rapidly decreasing sequences of numbers (each function
corresponds to its Fourier coefficients). On the other hand, this
transform can be considered as a $C^{\infty}$-version of the Fourier-Mukai
transform.
Roughly speaking, given a symplectic manifold $M$ with a Lagrangian tori
fibration, one introduces a natural complex structure on the dual
fibration $M\dl$. We say that $M\dl$ is mirror dual to $M$.
Then our transform produces a holomorphic vector bundle on $M\dl$
starting from a Lagrangian submanifold $L$ of $M$ transversal to all fibers and
a local system on $L$. We prove that the Dolbeault complex of this
holomorphic vector bundle is isomorphic to some modification of
the de Rham complex of the local system on $L$.
In the case of an elliptic curve, we check that all holomorphic vector bundles
on $M\dl$ are obtained in this way. Also we construct a quasi-isomorphism
of our modified de Rham complex with the complex
that computes morphisms in the Fukaya category between $L$ and some 
fixed Lagrangian submanifold (which corresponds to the trivial line
bundle on $M\dl$). One can construct a
similar quasi-isomorphism for arbitrary pair of Lagrangian submanifolds in 
$M$ (which
are transversal to all fibers). The most natural way to do it would
be to use tensor structures on our categories. The slight
problem is that we are really dealing with dg-categories rather than
with usual categories and the axiomatics of tensor dg-categories
does not seem to be understood well enough.
Hence, we restrict ourself to giving a brief sketch of how these
structures look in our case in Sections \ref{RemTen}, \ref{RemTenAbs}, and 
\ref{RemTenRel}. 
It seems that to compare Fukaya complex with our modified de Rham
complex in higher-dimensional case we need a generalization
of Morse theory for closed $1$-forms (cf. \cite{N}, \cite{Pa}) together
with a version of the result of Fukaya and Oh in \cite{FO} comparing
Witten complex with Floer complex.

The study of mirror symmetry via Lagrangian fibrations originates from
the conjecture of \cite{SYZ} that all mirror dual pairs of Calabi-Yau 
are equipped with dual special Lagrangian tori fibration.
The geometry of such fibrations and their compactifications is studied
in \cite{G1}, \cite{G2} and \cite{H}. In particular, the construction
of a complex structure on the dual fibration can be found in these papers.
On the other hand, K.~Fukaya explains in \cite{F-ab} how to construct
a complex structure (locally) on the moduli space of Lagrangian submanifolds
(equipped with rank $1$ local systems)
of a symplectic manifold $M$, where Lagrangian submanifolds are considered
up to Hamiltonian diffeomorphisms of $M$. Presumably these two constructions
are compatible and one can hope that for some class of Lagrangian submanifolds
the speciality condition picks a unique representative in each orbit of
Hamiltonian diffeomorphisms group. Our point of view is closer to that
of Fukaya: we do not equip our symplectic manifold with a complex structure,
so we cannot consider special geometry. 
However, we do not consider the problem of
compactifying the dual fibration and we do not know how to deal with
Lagrangian submanifolds which intersect some fibers non-transversally.
So it may well happen that special geometry will come up in relation
with one of these problems. 

The simplest higher-dimensional case in which our construction
can be applied is that of a (homogeneous) symplectic torus
equipped with a Lagrangian fibration by
\select{affine} Lagrangian subtori. The corresponding construction of the
mirror complex torus and of holomorphic bundles associated with
affine Lagrangian subtori intersecting fibers transversally
coincides with the one given by Fukaya in \cite{F-ab}.
However, even in this case the homological mirror conjecture
still seems to be far from reach (for dimensions greater than $2$).
Note that the construction of the mirror dual complex torus to a
given (homogeneous) symplectic torus $T$ requires a choice of a 
linear Lagrangian subtorus in $T$. For different choices we obtain
different complex tori. The homological mirror conjecture would
imply that the derived categories on all these complex tori
are equivalent (to be more precise, some of these categories should be
twisted by a class in $H^2(T\dl,\cal O^*)$). This is indeed the case
and follows from the main theorem of \cite{P}. The corresponding
equivalences are generalizations of the Fourier-Mukai transform. 

While we were preparing this paper, N.~C.~Leung and E.~Zaslow told us that
they invented the same construction of a holomorphic bundle coming
from a Lagrangian submanifold. 

\subsection{Organization}

Section 1 contains the basic definitions and a sketch of the results of this 
paper.

In Section 2, we deal with a single real torus. We define the Poincar\'e bundle that 
lives on the product of our torus
and the dual torus, and then use it to define a modified Fourier transform, 
which in this simple case is just the
correspondence between bundles with unitary connections on a torus and sky-scraper 
sheaves on the dual torus. 

Section 3 contains generalization of these results to families of tori. We describe 
the holomorphic sections of a vector bundle on the 
complex side in terms of rapidly decreasing
sections of some bundle corresponding to its ``Fourier transform'' (notice that 
not every holomorphic vector
bundle has this ``Fourier transform''). Here we also analyze the case of elliptic 
curve.

Section 5 is devoted to interpreting the Floer cohomologies in our terms 
(i.e., using the spaces of rapidly
decreasing sections of some bundles). This result is valid for elliptic curves 
only.  

\subsection{Notation}

We work in the category of real $C^\infty$-manifolds. The words ``a bundle on a 
manifold $X$'' mean a
(finite-dimensional) $C^\infty$-vector bundle over $\C$ on $X$. 
We usually identify a vector bundle with the corresponding sheaf of $C^\infty$-
sections. 
For a manifold $X$, $T_X\to X$ (resp. $T\dl_X\to X$) is the real
tangent (resp. cotangent) bundle, $\Vect(X)$ (resp. $\Omega^1(X)$) is the space of 
complex vector fields 
(resp. complex differential forms). If $X$ carries a complex structure,
$T^{0,1}_X\subset T_X\otimes\C$ stands for the subbundle of anti-holomorphic vector
fields. 
$\Diff(X)$ is the algebra of differential operators on $X$ with 
$C^\infty(X)\otimes\C$-coefficients.

Let $F$ be a vector bundle on a manifold $X$, $\nabla_F:F\to 
F\otimes\Omega^1(X)$ a connection.
We define the \select{curvature} $\curv\nabla_F\in\Omega^2(X)\otimes\End(F)$ of 
$\nabla_F$ by the usual formula:
$\langle(\curv\nabla_F),\tau_1\wedge\tau_2\rangle=
(\nabla_F)_{[\tau_1,\tau_2]}-[(\nabla_F)_{\tau_1},(\nabla_F)_{\tau_2}]$ for any 
$\tau_1,\tau_2\in\Vect(X)$. Here 
$\langle\al,\al\rangle$ stand for the
natural pairing $\bigwedge^2\Omega^1(X)\times\bigwedge^2\Vect(X)\to 
C^\infty(X)$ defined by 
$\langle\mu_1\wedge\mu_2,\tau_1\wedge\tau_2\rangle=
\langle\mu_1,\tau_2\rangle\langle\mu_2,\tau_1\rangle-
\langle\mu_1,\tau_1\rangle\langle\mu_2,\tau_2\rangle$.
 
A \select{local system} $\cL$ on a manifold $X$ is a vector bundle $F_\cL$ 
together with a 
connection $\nabla_\cL$ in $F_\cL$ such that $\curv\nabla_\cL=0$ (in other 
words, $\nabla_{\cL}$ is \select{flat}). 
The fiber $\cL_x$ of a local system $\cL$ over $x\in X$ equals the fiber 
$(F_\cL)_x$. 
For any $x\in X$, a local system $\cL$ defines the \select{monodromy} 
$\mon(\cL,x):\pi_1(X,x)\to GL(\cL_x)$.
We say that a local system is \select{unitary} if for any $x\in X$, there is
a Hermitian form on $\cL_x$ such that $\mon(\cL,x)(\gamma)$ are unitary
for all $\gamma\in\pi_1(X,x)$ (it is enough to check the condition for one
point on each connected component of $X$).

For a manifold $X$ and $\tau\in\Omega^1(X)$, we denote by $O_X(\tau)$ the 
trivial line bundle together 
with the 
connection $\nabla=d+\tau$. In particular, $O_X:=O_X(0)$ stands for the trivial 
local system on $X$.

\section{Main results}
\label{MainSec}

\subsection{}

Let $(M,\omega)$ be a symplectic manifold, $p:M\to B$ a surjective smooth map 
with Lagrangian fibers. 
Suppose that the fibers
of $M\to B$ are isomorphic to a torus $(\R/\Z)^n$. Fix a Lagrangian section 
$0_M:B\to M$. We call such collection 
$(p:M\to B,\omega,0_M)$ (or, less formally, the map $p:M\to B$) a 
\select{symplectic family of tori}.

The symplectic form induces a natural flat connection on $T_B$ 
(using the canonical isomorphism
$R^1p_*\R=T_B$) and an identification 
$M=T\dl_B/\Gamma$, where $\Gamma$ is a horizontal lattice in $T\dl_B$ ($\Gamma$ 
is dual to 
$\Gamma\dl:=R^1p_*\Z\subset R^1p_*\R=T_B$). This identification agrees with the 
symplectic structure, so $\Gamma\subset T\dl_B$ is Lagrangian.

Hence the connection on $T_B$ is symmetric (in 
the sense of \cite{Ml}). Recall that
a connection $\nabla$ on $T_B$ is called \select{symmetric} if
$\nabla_{\tau_1}(\tau_2)-\nabla_{\tau_2}(\tau_1)=
[\tau_1,\tau_2]$ for any $\tau_1,\tau_2\in\Vect(B)$. 

\begin{Rem} In particular, we see that $M\to B$ is locally (on $B$)
isomorphic to $(V/\Gamma)\times U$ for some vector space $V$, a lattice $\Gamma\subset V$,
and an open subset $U\subset V\dl$.
Besides, we see that the connection on $T_B$ induces a natural flat connection 
on $T_M$.
\end{Rem}

Consider the family of dual tori $M\dl:=T_B/\Gamma\dl$. The connection on $T_B$ 
yields a natural 
isomorphism $T_{M\dl}=(p\dl)^*T_B\oplus (p\dl)^*T_B$ such that the differential 
of $p\dl:M\dl\to B$ coincides with 
the first projection. So one can define a complex structure on $M\dl$ using the 
operator 
$J:T_{M\dl}\to T_{M\dl}:(\xi_1,\xi_2)\mapsto(-\xi_2,\xi_1)$. The complex 
manifold $M\dl$ is called 
the \select{mirror dual} of $M$.

For any torus $X=V/\Gamma$, the dual torus $X\dl=V\dl/\Gamma\dl$ can be 
interpreted as a moduli space of 
one-dimensional unitary local systems on $X$. So there is a natural universal 
$X\dl$-family of local systems
on $X$. We can interpret this family as a bundle with a connection on $X\times 
X\dl$ (see Section \ref{Poin} for details). 
If we apply this constructions to fibers of $p$, we get a canonical bundle $\cP$ 
on $M\times_BM\dl$ together with
a connection $\nabla_\cP$ on $\cP$ ($\nabla_\cP$ is not flat).

Suppose we are given a Lagrangian submanifold $i:L\hookrightarrow M$ which is 
transversal to fibers of $p$, and a local
system $\cL$ on $L$. We also assume that $p|_L:L\to B$ is proper. Define the 
\select{Fourier transform} of $(L,\cL)$ by 
the formula

\begin{equation}
\Four(L,\cL):=(p_{M\dl})_*
(((i\times id)^*\cP)\otimes((p_L)^*\cL))
\label{relFour}
\end{equation}

Here $p_{M\dl}:L\times_B M\dl\to M\dl$, $(i\times id):L\times_BM\dl\to 
M\times_BM\dl$, and 
$p_L:L\times_B M\dl\to L$ are the natural maps. The map $p_{M\dl}$ is a proper  
unramified covering, so
$\Four(L,\cL)$ is a bundle with connection on $M\dl$.

\begin{Th} (i) The $\overline\partial$-component of the connection on 
$\Four(L,\cL)$ is flat
(so $\Four(L,\cL)$ can be considered as a holomorphic vector bundle on $M\dl$);

(ii) If $B\simeq(\R/\Z)$, any holomorphic vector bundle on $M\dl$ is isomorphic 
to $\Four(L,\cL)$ for
some $(L,\cL)$. 
\label{dflat}
\end{Th}

\begin{Rem} There is an analogue of the above theorem for the case when
the fibration does not have a global Lagrangian section. In this case,
the dual complex manifold $M\dl$ carries a canonical cohomology class
$e\in H^2(M\dl,O_{M\dl}^*)$, hence one has the corresponding twisted category
of coherent sheaves (cf. \cite{G}). The analogue of $\Four(L,\cL)$ will
be an object in this twisted category. We will consider this generalization
in more details elsewhere. Also it would be interesting to find an
analogue of our construction for Lagrangian foliations. In the case of a
torus, this should lead to the functor considered by Fukaya in \cite{F-nc}.
\end{Rem}

\subsection{}

Let $(L,\cL)$ be as before.
 
Consider the natural map $u:T\dl_B\to M$ (the ``fiberwise universal cover''). Set 
$\tL:=u^{-1}(L)$. 
Denote by $u_L^*\cL$ the pull-back of $\cL$ to $\tL$ and by $\tau$ the restriction 
of the canonical $1$-form from
$T\dl_B$ to $\tL$. Since $\tL\subset T\dl_B$ is Lagrangian, $\tau$ is 
closed, so by adding
$-2\pi\tau$ to the connection on $u_L^*\cL$ we get a new local system 
$\tcL:=(u_L^*\cL)\otimes O_{\tL}(-2\pi\tau)$.

Denote by $C^\infty(\tcL)$ the space of $C^\infty$-sections of $\tcL$. 
Since $\tL\to B$ is an unramified covering,
we have an embedding $\Diff(B)\to\Diff(\tL)$. Set

\begin{equation}
\cS(\tcL):=\{s\in C^\infty(\tcL)|Ds\mbox{ is rapidly decreasing for any 
}D\in\Diff(B)\}
\end{equation}

Here a section $s$ of $\tcL$ is 
called \select{rapidly decreasing} if $\lim_{||g||\to\infty, 
g\in\Gamma_x}s((x,\tau+g))||g||^k=0$
for any $(x,\tau)\in L\times_M T\dl_B=\tL$ and $k>0$. Here $\Gamma_x$ stands for the 
fiber of $\Gamma\subset T\dl_B$ over $x\in B$. Since $s((x,\tau+g))\in 
\tcL_{(x,\tau+g)}=\cL_{(x,\tau)}$, the definition makes sense. Besides, it does not 
depend on the choice of a norm
$||\al||$ on $T\dl B$. Clearly, $\cS(\tcL)$ is a $\Diff(B)$-module.

\begin{Th} The de Rham complex $\DR(\tcL)$ of the $\Diff(B)$-module 
$\cS(\tcL)$ is isomorphic
to the Dolbeault complex of $\Four(L,\cL)$.
\label{ThdR}
\end{Th}

\subsection{}
\label{FukSec}

Suppose $B\simeq\R/\Z$. Fix an orientation on $B$.

Let $L,\cL$ be as before. Moreover, we suppose that $\cL$ is \select{quasi-unitary}, 
that is, for any $x\in L$ all eigenvalues of $\mon(\cL,x)$ are of absolute value $1$ 
(it follows from Lemma \ref{QU} that this condition
is not too restrictive). We also assume that $L$ meets the zero section 
$0_M(B)\subset M$
transversally.

As before, $\tL=u^{-1}(L)\subset T\dl_B$. Suppose $\tilde c\in\tL$ lies on the zero section 
$0_{T\dl_B}(B)\subset T\dl_B$. Then in a neighborhood of $\tilde c$, 
$\tL\subset T\dl_B$ is the graph of some $\mu\in\Omega^1(B)$, $d\mu=0$. Denote by $b\in B$ the image
of $\tilde c\in\tL$. In a neighborhood of $b$, $\mu=df$ for some $f\in C^\infty(B)$. We say that $\tilde c$
is \select{positive} (resp. \select{negative}) if $f$ has a local minimum (resp. maximum) at $b$.
Denote by $\{\tilde c_k^+\}\subset\tL$ (resp. $\{\tilde c_l^-\}\subset\tL$) the set of all positive (resp. negative) 
points of intersection with the zero section.

Let $\gamma\subset\tL$ be an arc with endpoints $\tilde c_k^+$ and $\tilde c_l^-$. We say that $\gamma$
is \select{simple} if it does not intersect the zero section.
Denote by $M(\gamma):\tcL_{\tilde c_k^+}\to\tcL_{\tilde c_l^-}$ the monodromy of $\tcL$ along 
$\gamma$ (the monodromy is the product of the monodromy of $u_L^*\cL$ and $\exp(2\pi A)$,
where $A$ is the oriented area of the domain restricted by $\gamma$ and the zero section).
Set $d(\gamma)=M(\gamma)$ if the direction from $\tilde c_k^+$ to $\tilde c_l^-$ along 
$\gamma$ agrees with the orientation of $B$, and $d(\gamma)=-M(\gamma)$ otherwise.

Set $F^0:=\oplus_k\tcL_{\tilde c_k^+}$, $F^1:=\oplus_l\tilde \cL_{\tilde c_l^-}$. Consider the operator 
$d:F^0\to F^1$ whose ``matrix elements'' are
$d_{kl}:\tcL_{\tilde c_k^+}\to\tcL_{\tilde c_l^-}=\sum_\gamma d(\gamma)$. Here the sum is taken 
over all simple arcs $\gamma$
with endpoints $\tilde c_k^+$, $\tilde c_l^-$ (there are at most two of them).
 
\begin{Rem} Since $L$ meets the fibers of $M\to B$ transversely, there is a 
canonical choice of lifting
$L\to M$ to $L\to\widetilde{GrL}(T_M)$. Here $GrL(T_M)\to M$ is the fibration whose fiber 
over $m\in M$ is the manifold
of Lagrangian subspaces in $T_M(m)$ (the \select{Lagrangian Grassmanian} of 
$T_M(m)$), $\widetilde{GrL}(T_M)\to GrL(T_M)$ is its fiberwise universal cover. This implies that the corresponding 
Floer cohomologies are
equipped with a natural $\Z$-grading.
Since $L$ is also transversal
to the zero section $0_M(B)\subset M$, we may compute the space (or, more precisely, the complex) of morphisms 
for the pair $L$, $0_M(B)$ in the Fukaya category. It is easy to see 
that the complex coincides with
$\cF(\cL):F^0\to F^1$. 
\end{Rem}

\begin{Th}
The complex $\cF(\cL)$ is quasi-isomorphic to $\DR(\tcL)$.
\label{ThFuk}
\end{Th} 

{\it Construction of a quasi-isomorphism $\cF(\cL)\to\DR(\tcL)$}.
Consider distributions with values in $\tcL$ that 
are rapidly decreasing smooth sections of $\tcL$ outside some compact set.
Let $\cS(\tcL)^D$ be the space of such distributions. 
Denote by $\DR(\tcL)^D$ the de Rham
complex associated with the $\Diff(B)$-module $\cS(\tcL)^D$. The inclusion 
$\cS(\tcL)\hookrightarrow\cS(\tcL)^D$
induces a quasi-isomorphism $\DR(\tcL)\to\DR(\tcL)^D$. Now let us define a 
morphism $\cF(\cL)\to\DR(\tcL)^D$.

For $\tilde c_k^+$, denote by $C_k^+$ the maximal (open) subinterval 
$I\subset\tL$ such that 
$\tilde c_k^+\in C_k^+$ and $\tilde c_l^-\notin C_k^+$ for any $l$ ($I$ may be 
infinite). 
The morphism $F^0\to\cS(\tcL)^D$
sends $v\in\tcL_{\tilde c_k^+}$ to $f$ such that $f$ vanishes outside $C_k^+$, 
$f$ is horizontal on $C_k^+$,
and $f(\tilde c_k^+)=v$. The morphism $F^1\to\cS(\tcL)^D\otimes\Omega^1(B)$ 
sends $v\in\tcL_{\tilde c_l^-}$ to
$v\otimes\delta_{\tilde c_l^-}$. Here $\delta_{\tilde c_l^-}$ is the delta-function at $\tilde c_l^-$.

\begin{Rem} All this machinery works in a more general situation. Namely, 
we can consider a
symplectic family of tori $M\to B$ together with a closed purely imaginary horizontal 
form $\omega^{I}$. Then we can work with the category of submanifolds 
$L\subset M$ together with a bundle $\cL$ on $L$ and a connection $\nabla_{\cL}$
such that $L\to B$ is a finite unramified covering and  
$\curv\nabla_\cL=2\pi(\omega+\omega^I)|_L$.
\end{Rem}

\subsection{}
\label{RemTen}

The pairs $(L,\cL)$ of the kind considered above form a category. One can define the (fiberwise) 
convolution product in this
category using the group structure on the fibers. However, the support of the 
convolution product does not
need to be a smooth Lagrangian submanifold, so to have a tensor category, we 
have to consider a slightly
different kind of objects (see Section \ref{ImLag}). 

After these precautions, we have a tensor category $\Sky(M/B)$. 
One easily sees that there is a canonical (i.e., functorial) choice
of the dual object $c\dl$ for any $c\in\Sky(M/B)$. For any $c\in\Sky(M/B)$, we 
have the de Rham complex 
$\DR(c)$ (defined in a way similar to what we do for $(L,\cL)$). Now we can use 
these data 
to define another ``category'' $\tSky(M/B)$:
we set $Ob(\tSky(M/B)):=Ob(\Sky(M/B))$, 
$\Hom_{\tSky(M/B)}(c_1,c_2):=\DR(c_2\star c_1\dl)$, where $\star$
stands for the convolution product. It is not a ``plain'' category, but a 
``dg-category''. Similarly,
the category of holomorphic vector bundles on $M\dl$ has a structure of a tensor 
dg-category 
(the morphism complex from $L_1$ to $L_2$ is the Dolbeault complex of 
$L_2\otimes L_1\dl$).  
Then the isomorphism of Theorem \ref{ThdR} induces a fully
faithful tensor functor between tensor dg-categories.

\section{Fourier transform on tori}

\subsection{Poincar\'e bundle}
\label{Poin}

Let $X$ be a torus (that is, a compact commutative real Lie group). 
Then $X=V/\Gamma$ for $V:=H_1(X,\R)$, $\Gamma:=H_1(X,\Z)$. 
The \select{dual torus} is $X\dl:=V\dl/\Gamma\dl$
($V\dl:=\Hom(V,\R)=H^1(X,\R)$, $\Gamma\dl:=\Hom(\Gamma,\Z)=H^1(X,\Z)$).

\begin{Df} A \select{Poincar\'e bundle} for $X$ is a line bundle $\cP$ on 
$X\times X\dl$
together with a connection $\nabla_\cP$ such that the following conditions are 
satisfied:

$(i)$ $\nabla_\cP$ is flat on $X\times\{x\dl\}$, and the monodromy is
$\pi_1(X)=H_1(X,\Z)\to U(1):
\gamma\mapsto \exp(2\pi \sqrt{-1}\langle x\dl,\gamma\rangle)$ 
(we denote by $\langle\al,\al\rangle$ not only the natural pairing 
$V\dl\times V\to\R$, but also the induced pairings $\Gamma\dl\times 
V/\Gamma\to\R/\Z$
and $V\dl/\Gamma\dl\times\Gamma\to\R/\Z$);
 
$(i\dl)$ $\nabla_\cP$ is flat on $\{x\}\times X\dl$, and the monodromy is
$\pi_1(X\dl)=H^1(X,\Z)\to U(1):\gamma\dl\mapsto \exp(-2\pi \sqrt{-1}
\langle \gamma\dl,x\rangle)$;

$(ii)$ For any $(x,x\dl)\in X\times X\dl $, $\delta v\in V=T_xX$, 
$\delta v\dl\in V\dl=T_{x\dl} X\dl$, we have
$\langle\curv(\nabla_\cP),\delta v\wedge \delta v\dl\rangle=-2\pi\sqrt{-1}\langle\delta v\dl,\delta v\rangle$.     
\end{Df}

Clearly, $(\cP,\nabla_\cP)$ is defined up to an isomorphism by $(i)$, $(i\dl)$, 
$(ii)$. Furthermore, we always fix an identification $\iota:\cP_{(0,0)}\iso\C$, 
so the collection $(\cP,\nabla_\cP,\iota)$ is defined up to a canonical 
isomorphism.

A Poincar\'e bundle allows us to identify $X\dl$ with the moduli space of 
unitary
local system on $X$ (and vice versa).

\begin{Rem} Suppose $V$ carries a complex structure $J:V\to V$. Define the complex
structure on $V\dl$ using $-J\dl$. Then $X$, $X\dl$, and $X\times X\dl$ are complex manifolds.
Let $\cP$ be a Poincar\'e bundle for $X$. It is easy to see that $\nabla_\cP$ is 
``flat in $\overline\partial$-direction'' (i.e., the $\curv\nabla_{\cP}$ vanishes on 
$\bigwedge^2T^{0,1}_{X\times X\dl}$). Hence $\cP$ can be considered as a holomorphic line
bundle on $X\times X\dl$. Actually, $\cP$ is in this case isomorphic to the ``complex'' Poincar\'e
bundle (i.e., the universal bundle that comes from the interpretation of $X\dl$ as a moduli space
of holomorphic line bundles on $X$).
\end{Rem} 

The following lemma is straightforward.

\begin{Lm} Consider the local system 
$F:=O_{V\times X\dl}(2\pi\sqrt{-1}\langle dx\dl,v\rangle)$. Here $dx\dl\in\Omega^1(X\dl)\otimes V\dl$
is the natural form with values in $V\dl$.
Lift the natural action of $\Gamma=H_1(X,\Z)$ on $V$ to $F$ by 
$(g(f))(v,x\dl)=\exp(-2\pi\sqrt{-1}\langle x\dl,g\rangle)f(v-g,x\dl)$. 
Then the corresponding line bundle with connection on $X\times X\dl$
is a Poincar\'e bundle. 
%\QED
\end{Lm}

Consider the natural projection 
$u\times id:V\times X\dl\to X\times X\dl$. 
Then $(u\times id)^*\cP$
is identified with $F$. 
We denote by $\ex{v}{x\dl}$ the section of $(u\times id)^*\cP$ that corresponds 
to
$1\in F$.

\begin{Rem}
Let $\cP$ be a Poincar\'e bundle for $X$, $\sigma':X\dl\times X\to
X\times X\dl:(x\dl,x)\mapsto(-x,x\dl)$. Then $(\sigma')^*\cP$ is a Poincar\'e
bundle for $X\dl$. 
\end{Rem}

\subsection{Sky-scraper sheaves}

Given a finite set $S\subset X$  and (finite-dimensional) $\C$-vector spaces 
$F_s$ for all $s\in S$,
we can define the corresponding \select{(finite semisimple) sky-scraper sheaf} 
$F$ on $X$ by
$F(U)=\oplus_{s\in S\cap U}F_s$ for $U\subset X$. Denote by $\Sky(X)$ the 
category of sky-scraper sheaves on $X$
($\Sky(X)$ is a full subcategory of the category of sheaves of vector spaces on 
$X$). 
Any sky-scraper sheaf is naturally a 
$C^\infty(X)$-module, and morphisms of sky-scraper sheaves agree with the 
action of $C^\infty(X)$.

For $F\in\Sky(X)$, define the Fourier transform of $F$ by

\begin{equation}
\Four F:=(p_{X\dl})_*((p_X^* F)\otimes\cP)
\end{equation} 

Here $p_X:X\times X\dl\to X$ and $p_{X\dl}:X\times X\dl\to X\dl$ are
the natural projections.

$\Four F$ is a locally free sheaf of rank $\dim H^0(X,F)$, so we interpret 
$\Four F$ as a vector bundle
on $X\dl$. The connection $\nabla$ on $\cP$ induces
a flat unitary connection on $\Four F$. So $\Four$ can be considered as a 
functor 
$\Sky(X)\to\LU(X\dl)$, where $\LU(X\dl)$ is the category of unitary local 
systems on $X\dl$. 
This functor is an equivalence of categories.

\subsection{Rapidly decreasing sections}

For a sheaf $F\in\Sky(X)$, set $\tF:=u^*F$, where $u:V\to X$ is the universal 
cover. 
The group $\Gamma:=H_1(X,\Z)$ 
acts on $V=H_1(X,\R)$ and $\tF$ is $\Gamma$-equivariant. We say that a section 
$s\in H^0(V,\tF)$ is 
rapidly decreasing if $\lim_{||g||\to\infty,g\in\Gamma} s(x+g)||g||^{k}=0$ for any 
$x\in V$, $k>0$ (the definition does not depend 
on the choice of a norm $||\al||$ on $V$). Denote
by $\cS(\tF)$ the space of all rapidly decreasing sections of $\tF$.

Take $F\in\Sky(X)$, $f\in\cS(\tF)$. Set

\begin{equation}
\Four_F f(x\dl)=\sum_{v\in V} f(v)\ex{v}{x\dl}
\label{sectFour}
\end{equation}

The following lemma is clear:

\begin{Lm} Let $F\in\Sky(X)$. Then $\Four_F:\cS(\tF)\to C^\infty(\Four(F))$ 
is an isomorphism. Here $C^\infty(\Four(F))$
is the space of $C^\infty$-sections of the local system $\Four(F)$.
%\QED
\label{cSLm}
\end{Lm}

\subsection{Convolution}
\label{RemTenAbs}

For $F_1,F_2\in\Sky(X)$, one can define their \select{convolution product}
by $F_1\star F_2:=m_*((p_1^*F_1)\otimes(p_2^*F_2))$, where $m,p_1,p_2:X\times 
X\to X$ are the group law, 
the first projection, and the second projection respectively. This gives a 
structure of a tensor category on 
$\Sky(X)$ (the unit, dual element, and commutativity and associativity 
isomorphisms are easily defined). 
Then $\Four:\Sky(X)\to\LU(X\dl)$ is naturally a tensor functor (the tensor 
structure on $\LU(X\dl)$ is the ``usual''
tensor product). Moreover, $\Four(F_1\star F_2)=\Four(F_1)\otimes\Four(F_2)$ for any $F_1,F_2\in\Sky(X)$.

Besides, it is easy to define the natural convolution product 
$\cS(\widetilde{\star\vphantom{F}}):\cS(\tF_1)\otimes\cS(\tF_2)\to\cS(\widetilde{(F_1\star 
F_2)})$.
This makes $\cS(\widetilde{\vphantom{F}\al})$ a tensor functor.
One can check that $\Four_\al:\cS(\widetilde{\vphantom{F}\al})
\to C^\infty(\Four(\al))$ is actually an isomorphism
of tensor functors (i.e., for any $F_1,F_2\in\Sky(X)$ the diagram

\begin{equation}
\label{comm}
\begin{array}{ccc}
\cS(\tF_1)\otimes\cS(\tF_1)&\iso&C^\infty(\Four(F_1))\otimes 
C^\infty(\Four(F_2))\\
\downarrow&\;&\downarrow\\
\cS(\widetilde{F_1\star F_2})&\iso&C^\infty(\Four(F_1)\otimes\Four(F_2))
\end{array}
\end{equation}

commutes).

\begin{Exm} Let $F$ be the \select{unit object} in $\Sky(X)$ (i.e., $\supp 
F=\{0\}$ and $F_0=\C$). Then 
$\tF$ is a trivial sheaf on $\Gamma=H^1(X,\Z)$. Clearly, $\Four(F)=O_X$ is the 
trivial local system on $X$. 
In this case, the isomorphism $(\Four_F)^{-1}:C^\infty(X\dl)\to\cS(\tF)$ maps 
any $C^\infty$-function to its Fourier 
coefficients. Since $F\star F=F$, the commutativity of (\ref{comm}) 
in this case is the 
well-known formula for the Fourier coefficients of the product.
\end{Exm}

\section{Relative sky-scraper sheaves}

Let $p:M\to B$ be a symplectic family of tori. 
In this section, we construct ``relative versions'' of the objects from the 
previous section.

\subsection{}

\label{ImLag}

A \select{transversally immersed Lagrangian manifold} is a couple $(L,i)$, 
where $i:L\to M$ is a morphism of $C^\infty$-manifolds such that
$p\circ i:L\to B$ is a proper finite unramified covering and $i^*(\omega)=0$.

Consider the category $\Sky(M/B)$, whose objects are
triples $(L,i,\cL)$, where $(L,i)$ is a transversally immersed 
Lagrangian submanifold, and $\cL$ is a local system on $L$.

\begin{Rem}%???
Take any $(L_1,i_1,\cL_1),(L_2,i_2,\cL_2)\in\Sky(M/B)$. Consider $L_{1\to2}':=
L_1\times_M L_2$. Denote by $L_{1\to2}\subset L_{1\to2}'$ 
the maximal closed submanifold whose images in $L_1$, $L_2$ are open
(if $L_1$ and $L_2$ are just ``usual'' Lagrangian submanifolds, $L_{1\to2}$ is the 
union of common
connected components of $L_1$ and $L_2$). Let $p_1:L_{1\to2}\to L_1$, 
$p_2:L_{1\to2}\to L_2$ be the 
natural projections.
By definition, morphisms from $(L_1,i_1,\cL_1)$ to $(L_2,i_2,\cL_2)$ are 
horizontal morphisms 
$p_1^*\cL_1\to p_2^*\cL_2$. The composition is defined in the natural way. 
\end{Rem}

\subsection{}

Let $p\dl:M\dl\to B$ be the mirror dual of $M\to B$. Take 
$(L,i,\cL)\in\Sky(M/B)$. 
One can easily define the (relative) Poincar\'e bundle $\cP$ on $M\times_BM\dl$. 
It carries a natural
connection $\nabla_\cP$. Now define the (fiberwise) Fourier transform 
$\Four(L,i,\cL)$ by the formula (\ref{relFour}).

\begin{proof}[Proof of Theorem \ref{dflat}(i)]
The natural map
$L\times_BM\dl\to M\dl$ is an unramified covering, so the complex structure on $M\dl$
induces a complex structure on $L\times_BM\dl$. Let $(\cP_M,\nabla_{\cP_M})$ be the 
Poincar\'e bundle on $M$. It is enough to prove $\curv(\nabla_{\cP_M})$ vanishes on 
$T_{L\times_BM\dl}^{0,1}$. 

The statement is local on $B$, so we may assume $M=X\times B$, $M\dl=X\dl\times B$ for
a torus $X$. Denote by $p_{X\times X\dl}:M\times_BM\dl\to X\times X\dl$ the natural projection,
and by $(\cP_X,\nabla_{\cP_X})$ the Poincar\'e bundle of $X$.  
$p_{X\times X\dl}^*(\cP_X,\nabla_{\cP_X})=(\cP_M,\nabla_{\cP_M})$, so 
$\curv\nabla_{\cP_M}=p_{X\times X\dl}^*\curv\nabla_{\cP_X}$. Since $\curv\nabla_{\cP_X}$ is a scalar multiple 
of the natural symplectic form on $X\times X\dl$, it is enough to notice that
$p_{X\times X\dl}$ maps $T^{0,1}_{L\times_BM\dl}(x)$ to a Lagrangian subspace of
$T_{X\times X\dl}(p_{X\times X\dl}(x))\otimes\C$ for any $x\in L\times_BM\dl$.
\end{proof}

\subsection{Proof of Theorem \ref{ThdR}}

Consider the ``fiberwise universal cover'' $u:T\dl_B\to M$. For any 
$(L,i,\cL)\in\Sky(M/B)$, set
$\tL:=L\times_MT\dl_B$. 
Recall that $\tcL=u_L^*(\cL)\otimes O_\tL(-2\pi\tau)$, where $u_L:\tL\to L$ is 
the natural homomorphism and $\tau$
is the pull-back of the natural $1$-form on $T\dl_B$. 

For any $D\in\Diff(B)$, we consider its pull-back
$\tilde p^*D\in\Diff(\tL)$ (since $\tilde p:\tL\to B$ is an unramified covering, the pull-back 
is well defined). 
Since $\tcL$ carries a canonical
flat connection, we can apply $\tilde p^*D$ to $s\in C^\infty(\tcL)$. 
Denote by $\cS(\tcL)$ the set of all sections $s\in C^\infty(\tcL)$ such that
$(\tilde p^*D)s$ is (fiberwise) rapidly decreasing for any $D\in\Diff(B)$. 
$\cS(\tcL)$ is a $\Diff(B)$-module.

Just like in the ``absolute'' case (Lemma \ref{cSLm}), the Fourier transform 
(formula (\ref{sectFour})) yields a canonical isomorphism
$\cS(\tcL)\iso C^\infty(\Four(L,i,\cL))$ for any $L,i,\cL$.

The natural morphism
$(dp\dl)\otimes\C:T_{M\dl}\otimes\C\to(p\dl)^*T_B\otimes\C$ induces an 
isomorphism 
$T^{0,1}_{M\dl}\iso(p\dl)^*T_B\otimes\C$. 
So we have an embedding of Lie algebras 
$\Vect(B)\to C^\infty(T^{0,1}_{M\dl})\subset\Vect(M\dl)$.
The Lie algebra $C^\infty(T^{0,1}_{M\dl})$ of anti-holomorphic vector fields 
acts on 
$C^\infty(\Four(L,i,\cL))$ (by Theorem \ref{dflat}(i)),
so $C^\infty(\Four(L,i,\cL))$ has a natural structure of a $\Diff(B)$-module. 
One easily checks that the de Rham
complex associated with this $\Diff(B)$-module is identified with the Dolbeault 
complex of $\Four(L,i,\cL)$.

The following lemma implies Theorem \ref{ThdR}.

\begin{Lm} The isomorphism $\cS(\tL)\iso C^\infty(\Four(L,i,\cL))$ agrees with 
the structures of $\Diff(B)$-modules.
\end{Lm}
\begin{proof}
Again, we may assume $M=B\times X$ for a torus $X=V/\Gamma$. Consider the natural maps
$p_{V\times X\dl}:\tL\times_BM\dl\to T\dl_B\times_B M\dl\to V\times X\dl$,
$p_M:\tL\times_BM\dl\to T\dl_B\to M$, and $p_{T\dl_B}:\tL\times_BM\dl\to
T\dl_B\times_BM\dl\to T\dl_B$. 

$\ex{v}{x\dl}$ can be considered as  a horizontal section of $p_M^*(\cP_M)\otimes 
p^*_{V\times X\dl}(O_{V\times X\dl}(-2\pi\sqrt{-1}\langle dx\dl,v\rangle))$. Now the statement follows
from the fact that $1$ is a holomorphic section of 
$O_{\tL\times_BM\dl}(-2\pi p_{T\dl_B}^*\tau-2\pi\sqrt{-1}p_{V\times X\dl}^*\langle dx\dl,v\rangle)$ 
(i.e., the $\overline\partial$ component of the connection vanishes on $1$). Here $\tau$ stands for the
natural $1$-form on $T\dl_B$, and the complex structure on $\tL\times_BM\dl$ is that induced by 
$\tL\times_BM\dl\to M\dl$.
\end{proof}

\subsection{Proof of Theorem \ref{dflat}(ii)}

This result is actually proved in \cite{PZ}. Our proof is slightly different in 
that it makes use of connections.

Let $F$ be a holomorphic bundle on the elliptic curve $M\dl$. 
It is enough to consider the case of indecomposable $F$. 

The following statement is a reformulation of \cite[Proposition 1]{PZ} (which in 
turn is a consequence of
M.~Atiyah's results \cite{At}).

\begin{Pp} An indecomposable bundle $F$ on $M\dl$ is isomorphic to 
$\pi_{r,*}(L\otimes N)$,
where $\pi_r:M_r\dl\to M\dl$ is the isogeny corresponding to an (unramified) 
cover $B_r\to B$,
$L$ is a line bundle on $M\dl_r$, and $N$ is a unipotent bundle on $M\dl_r$ 
(i.e., $N$ admits a filtration
with trivial factors).
\end{Pp}

$\Four$ agrees with passing to unramified 
covers $B_r\to B$, besides, $\Four$ transforms the convolution product in $\Sky(M/B)$ to the
tensor product of holomorphic vector bundles (see Section \ref{RemTenRel} for the definition of the
convolution product). So it suffices to
consider the following cases:

{\it Case 1.} Let $F=l$ be a line bundle on $M\dl$. Our statement in this case 
follows from the following easy 
lemma:

\begin{Lm} $l$ carries a $C^\infty$-connection $\nabla_l$ such that the 
following conditions are satisfied:

$i)$ $\nabla_l$ agrees with the holomorphic structure on $l$ (i.e., the 
$\overline\partial$-component of 
$\nabla_l$ coincide with the canonical $\overline\partial$-differential);

$ii)$ The curvature $\curv\nabla_l$ is a horizontal (1-1)-form on $M\dl$ (in 
terms of the canonical connection);

$iii)$ The monodromies of $\nabla_l$ along the fibers of $M\dl\to B$  are 
unitary.
\end{Lm} 

{\it Case 2.} Let $F=N$ be a unipotent bundle on $M\dl$. To complete the proof, 
it is enough to notice that $N$ carries
a flat connection $\nabla_N$ such that $\nabla_N$ agrees with the holomorphic 
structure and
$\nabla_N$ is trivial along the fibers of $M\dl\to B$.

\subsection{Remarks on tensor dg-categories}
\label{RemTenRel}

For any $(L_1,i_1,\cL_1),(L_2,i_2,\cL_2)\in\Sky(M/B)$, set $L:=L_1\times_B L_2$, 
$\cL:=p_1^*(\cL_1)\otimes p_2^*(\cL_2)$ (here $p_i:L\to L_i$ is the natural 
projection). 
Consider the composition $i:=m\circ(i_1\times i_2):L_1\times_B L_2\to M\times_B 
M\to M$,
where $m:M\times_B M\to M$ is the group law $(x_1,x_2)\mapsto x_1+x_2$. 
Clearly, $(L,i,\cL)\in\Sky(M/B)$. 
$(L_1,i_1,\cL_1)\star(L_2,i_2,\cL_2):=(L,i,\cL)$ is the convolution product
of $(L_1,i_1,\cL_1)$ and $(L_2,i_2,\cL_2)$. The convolution product naturally 
extends to a structure of tensor category
on $\Sky(M/B)$ (the unit object, dual objects, and associativity/commutativity 
constraints are defined in a natural
way). Notice that there is a functorial choice of dual object. 

Just as in Section \ref{RemTenAbs}, the convolution product induces a functorial 
morphism of $\Diff(B)$-modules
$\cS(\tcL_1)\otimes\cS(\tcL_2)\to\cS(\tcL_3)$ for any 
$(L_1,i_1,\cL_1),(L_2,i_2,\cL_2)\in\Sky(M/B)$,
$(L_3,i_3,\cL_3):=(L_1,i_1,\cL_1)\star(L_2,i_2,\cL_2)$. So 
$\cS(\widetilde{\vphantom{L}\al})$ is a tensor 
functor from $\Sky(M/B)$ to the category of $\Diff(B)$-modules.

Just as we say in Section \ref{RemTen}, we define a tensor dg-category 
$\tSky(M/B)$ by setting
$Ob(\tSky(M/B)):=Ob(\Sky(M/B))$, $\Hom_{\tSky(M/B)}(c_1,c_2):=\DR(c_2\star 
c_1\dl)$.

\section{Connection with the Fukaya category}

\subsection{Hamiltonian diffeomorphisms}

In this section, we prove some results about tensor dg-category $\tSky(M/B)$. We 
do not use these facts anywhere,
so the part may be skipped. However, the results give some clarification to the 
connection between
$\tSky(M/B)$ and the original category considered by Fukaya \cite{Fuk}.

Fix $\mu\in\Omega^1(B)$ such that $d\mu=0$. $\mu$ can be considered as a section 
of $T\dl_B$. Denote by
$i_\mu:B\to M$ the image of this section via the fiberwise universal cover 
$T\dl_B\to M$. 
Set $c_\mu:=(B,i_\mu,O_B(2\pi\mu))\in\tSky(M/B)$. The following 
statement follows from the
definitions:

\begin{Pp} 
$c_\mu\simeq 1_{\tSky(M/B)}$ in $\tSky(M/B)$.
%\QED
\end{Pp}

Now let $A:M\to M$ be any symplectic diffeomorphism that preserves the fibration 
$M\to B$. It is easy to see
that $A$ preserves the action of $T\dl_B$ on $M$, so $A$ corresponds to some 
$\mu\in\Omega^1(B)$. Since
$A$ preserves the symplectic structure, $d\mu=0$. Now we can consider the 
``automorphism''  $c\mapsto c_\mu\star c$
of $\tSky(M/B)$. Note that if $(L',i_{L'},\cL')=(L,i_L,\cL)\star c_\mu$, then 
$i_{L'}L'=A(i_LL)$.

In particular, if $A$ is Hamiltonian (that is, there is $f\in C^\infty(B)$ such 
that $\mu=df$), 
we get the following statement:

\begin{Co}
The map $(L,i_L,\cL)\mapsto(L,A\circ i_L,\cL)$ extends to an automorphism of 
$\tSky(M/B)$.
\end{Co} 
\begin{proof}
It is enough to notice that $O_B(2\pi\mu)$ is a trivial local system 
if $\mu=df$,
so $c_\mu\simeq(B,i_\mu,O_B)$ and $(L,i_L,\cL)\star(B,i_\mu,O_B)=(L,A\circ 
i_L,\cL)$ for any
$(L,i_L,\cL)\in\tSky(M/B)$.
\end{proof}

From now on, we suppose that $B$ is a torus.

Denote by $\tSky(M/B)^{QU}$ the full subcategory of $\tSky(M/B)$ formed by 
triples 
$(L,i,\cL)$ with quasi-unitary $\cL$ (that is, all the eigenvalues of all the 
monodromy operators 
are of absolute value $1$).

\begin{Lm}
The natural embedding $\tSky(M/B)^{QU}\to\tSky(M/B)$ is an equivalence of %??? dg- 
categories.
\label{QU}
\end{Lm}
\begin{proof}
We should prove that for any $(L,i,\cL)\in\tSky(M/B)$ there is 
$(L',i',\cL')\in\tSky^{QU}(M/B)$ such that 
$(L,i,\cL)\simeq(L',i',\cL')$ in $\tSky(M/B)$. It is enough to prove this 
statement for indecomposable $\cL$ and 
connected $L$.

Choose a point $x\in L$. For $\gamma\in\pi_1(L)$, 
we denote the monodromy along $\gamma$ by $\mon(\gamma)\in GL(L_x)$.
For any loop $\gamma\in\pi_1(L)$, all the eigenvalues of $\mon(\gamma)$ are of 
the same absolute value
(otherwise $\cL$ is decomposable). 

Consider the homomorphism $\mu:\pi_1(L)\to\R_+:=\{a\in 
R|a>0\}:\gamma\mapsto|\det(\mon(\gamma))|^{-1/d}$. 
Since $L\to B$
is a finite covering, $\pi_1(L)\subset H_1(B,\Z)\subset H_1(B,\R)$ is a lattice. 
So $\mu$ induces $\log\mu\in\Hom(\pi_1(L),\R)=H^1(B,\R)$. 

Choose an invariant $1$-form $\tilde\mu$ on $B$ that represents $-\frac{\log\mu}{2\pi}\in H^1(B,\R)$. Clearly,
$(L,i,\cL)\star c_{\tilde\mu}\in\tSky(M/B)^{QU}$. 
\end{proof}

\begin{Rem} Suppose $M$ and $B$ are tori (in particular, they have a Lie group structure), and $p:M\to B$
is a group homomorphism.  
Assume also $\omega$ is translation invariant.
Let $i:L\to M$ be a transversally immersed Lagrangian submanifold. We say that $(L,i)$ is \select{linear}
if for any connected component $L_j\subset L$, one has $i(L_j)=m+L'$ for some $m\in M$ and some Lie subgroup 
$L'\subset M$. 
Consider the full subcategory $\tSky(M/B)^{LN}\subset\tSky(M/B)$ that consists of $(L,i,\cL)$ such that
$(L,i)$ is linear and $\cL$ is quasi-unitary. It can be proved (in a way similar to the proof of Lemma \ref{QU})
that $\tSky(M/B)^{LN}\to\tSky(M/B)$ is an equivalence of categories. 
\end{Rem}

\subsection{Proof of Theorem \ref{ThFuk}}

In this section we give a different construction of the quasi-isomorphism 
between $\cF(\cL)$ and $\DR(\tcL)$.
Identification of this quasi-isomorphism with that constructed in Section 
\ref{FukSec} is left to reader. 

Consider the de Rham complex 
$\DR(\cL):=\cS(\tcL)\toup{d}\cS(\tcL)\otimes_{C^\infty(B)}\Omega^1(B)$. 
Recall that $\{\tilde c_l^-\}\subset\tL$ is the set of all ``negative''
points whose images lie on the zero section $0_{T\dl_B}(B)\subset T\dl_B$. 

Denote by $F^{(0)}$ the set of all $f\in\cS(\tcL)$ such that $f$ is 
horizontal in a neighborhood of 
$\{\tilde c_l^-\}$ and by $F^{(1)}$ the set of all 
$\mu\in\cS(\tcL)\otimes_{C^\infty(B)}\Omega^1(B)$ such that
$\mu$ vanishes in a neighborhood of $\{\tilde c_l^-\}$. Since 
$d(F^{(0)})\subset F^{(1)}$,
we have a complex $\DR'(\cL):F^{(0)}\toup{d}F^{(1)}$. Moreover, the natural 
map $\DR'(\cL)\to\DR(\cL)$ is a 
quasi-isomorphism.

Now let $\tL_j\subset\tL$ be a connected component of $\tL\setminus\{\tilde 
c_l^-\}$. Set $\cS(\tcL)_j:=
\{s\in\cS(\tcL):\supp s\subset\tL_i\}$, 
$F^{(1)}_j:=\cS(\tcL)_j\otimes_{C^\infty(B)}\Omega^1(B)$. Denote by
$F_j'$ the set of all sections $s\in C^\infty(\tcL|_{\tL_j})$ such that 
$\supp s$ is contained in a compact set
$C\subset L$ and $s$ is horizontal in some neighborhood of $\{\tilde x_l^-\}$. Set
$F^{(0)}_j:=\cS(\tcL)_j+F_j'$. Clearly $d$ yields a morphism $d_j: F 
_j^{(0)}\to F _j^{(1)}$, so we have
a complex $\DR(\cL)_j: F_j^{(0)}\toup{d_j} F_j^{(1)}$. 

The restriction map induces a morphism of complexes $\DR'(\cL)\to\DR(\cL)_j$. So 
we have a map
$\DR'(\cL)\to\oplus_j\DR(\cL)_j$. Moreover, one can see that this map is 
included into a short exact sequence

\begin{equation}
0\to \DR'(\cL)\to\oplus_j\DR(\cL)_j\to F^1\to 0 
\label{exseq}
\end{equation}

(see Section \ref{FukSec} for the definition of $F^1$).

\begin{Pp} $i)$ The map $d_j: F_j^{(0)}\to F_j^{(1)}$ is surjective for any $j$;

          $iia)$ If the image of $\tL_j\subset\tL$ does not intersect $0_{T\dl_B}(B)\subset T\dl_B$, 
                 then $d_j$ is injective;

          $iib)$ Suppose the image of $\tilde c_k^+\in\tL_j$ lies on $0_{T\dl_B}(B)\subset T\dl_B$.
                 Set $F_j'':=\ker d_j$. Then the map 
$F''_j\to\tcL_{\tilde c_k^+}:s\mapsto s(\tilde c_k^+)$ is bijective.
\end{Pp}

\begin{proof}

Choose an isomorphism $t:B\iso\R/\Z$. We may assume that $t$ agrees with the 
natural connection on $T_B$.
There are three possibilities:

\select{Case 1:} $M\to B\toup{t}\R/\Z$ induces an isomorphism 
$t:\tL_j\iso\R/m\Z$ for some m.

In this case, the image of $\tL_j$ does not intersect $0_{T\dl_B}(B)\subset T\dl_B$.
It is easy to see that 
the monodromy of $\tcL|_{\tL_j}$ does not have $1$ as its eigenvalue, hence 
the de Rham cohomology groups
of $\tcL|_{\tL_j}$ vanish. So $d_j$ is bijective and $i)$, $iia)$ follow.

\select{Case 2:} $M\to B\toup{t}\R/\Z$ induces an isomorphism $t:\tL_j\iso 
(t_1,t_2):=\{t\in\R:t_1<t<t_2\}$ for some
$t_1,t_2\in\R$.

In this case, there is a unique $\tilde c_k^+\in\tL_j$ whose image lies on $0_{T\dl_B}\subset T\dl_B$.
Besides, 
$F_j^{(1)}=H^0_c(\tL_j,\tcL|_{\tL_j}\otimes_{C^\infty(\tL_j)}\Omega^1(\tL_j))$ 
(here $H^0_c$ stands for the space of sections with
compact support). Now $i)$, $iib)$ are obvious.

\select{Case 3:} $M\to B\toup{t}\R/\Z$ induces an isomorphism $t:\tL_j\iso 
(t_1,t_2)$ where either $t_1=-\infty$,
or $t_2=\infty$ (or both). Without loss of generality, we assume $t_1=-\infty$.

Denote by $\tau\in\Omega^1(\tL_j)$ the pull-back of the natural $1$-form on $T_B$. 
It is easy to see there are (unique) $a,b\in\R$ ($a\ne0$) such that  
$\tau_0:=-2\pi\tau-(ta+b)dt$ is ``bounded''
in the following sense: there is $C\in\R$ such that for any connected closed 
subset $U\subset\tL_j$ we have

\begin{equation}
\left|\int_U\tau_0\right|<C
\label{eqbound}
\end{equation}

Choose an isomorphism $\phi:\tcL|_{\cL_j}\iso(O_{\cL_j}((at+b)dt))^n$ 
(where $n$ is the dimension of 
$\tcL|_{\cL_j}$). Set $\widehat{\frac{d}{dt}}:=\frac{d}{dt}+at+b$. Denote by 
$\hat\cS(t_1,t_2)$ the space of 
$f\in C^\infty(t_1,t_2)$ such that $\lim_{t\to\infty}x^l\frac{d^kf}{dt^k}=0$ 
for any $k,l\ge0$. 
If $t_2<\infty$, we 
denote by $\hat\cS^0(t_1,t_2)\subset\hat\cS(t_1,t_2)$ (resp.  
$\hat\cS^1(t_1,t_2)\subset\hat\cS(t_1,t_2)$) the subspace
of functions $f$ such that $\widehat{\frac{d}{dt}}f=0$ (resp. $f=0$) in a neighborhood of
$t_2$. If $t_2=\infty$, we set $\hat\cS^0(t_1,t_2)=\hat\cS^1(t_1,t_2)=\hat\cS(t_1,t_2)$.
(\ref{eqbound}) implies that $\phi$
induces isomorphisms $F^{(0)}_j\iso(\hat\cS^0(t_1,t_2))^n$, 
$F^{(1)}_j\iso(\hat\cS^1(t_1,t_2))^ndt$.
The differential $d_j$ corresponds to $\widehat{\frac{d}{dt}}dt$.

There are two possibilities:

\select{Case 3a:} $a>0$, the image of $\tL_j$ intersects $0_{T\dl_B}(B)\subset T\dl_B$ in exactly one 
point. Without loss of generality we 
may assume this point corresponds to $t=0$. Now for any 
$g\in(\hat\cS^1(t_1,t_2))^n$, a generic
solution to $\widehat{\frac{d}{dt}}f=g$ is given by

\begin{equation}
f(x)=\exp(-(ax^2/2+bx))(\int_0^x g(t)\exp(at^2/2+bt)dt+C)
\end{equation}

where $C\in\C^n$. It is easy to see $f\in(\hat\cS^0(t_1,t_2))^n$ for any $C$.
$i)$ and $iib)$ follow.

\select{Case 3b:} $a<0$, the image of $\tL_j$ does not meet the zero section, 
$t_2<\infty$. For any
$g\in(\hat\cS^1(t_1,t_2))^n$ the formula

\begin{equation}
f(x)=\exp(-(ax^2/2+bx))(\int_{-\infty}^x g(t)\exp(at^2/2+bt)dt+C)
\end{equation}

gives a generic solution to $\widehat{\frac{d}{dt}}f=g$ ($C\in\C^n$).
 However, $f\in(\hat\cS^0(t_1,t_2))^n$ if and only if $C=0$. This implies
$i)$ and $iia)$.
\end{proof}

Hence $H^i(\oplus_j\DR(\cL)_j)=\begin{cases}F^0,j=0\cr 0,\mbox{ 
otherwise}\end{cases}$.
To complete the proof, it is enough to notice that the map 
$F^0=H^0(\oplus_j\DR(\cL)_j)\to F^1$ induced
by (\ref{exseq}) coincides with that defined in Section \ref{FukSec}.

\end{document}